\documentclass[a4paper,reqno]{amsart}
\usepackage{amssymb, amsmath, amscd}

\def\Pspec{\operatorname{Spec}}
\def\Tr{\operatorname{Tr}}
\def\vol{\operatorname{Vol}}

\newtheorem{thm}{Theorem}

\newtheorem{lem}[thm]{Lemma}

\newtheorem{defn}{Definition}

\makeatother
\begin{document}
\title[ Spectral geometry of eta-Einstein Sasakian manifolds]
{Spectral geometry of eta-Einstein \\ Sasakian manifolds}
\author{JeongHyeong Park}
\address{Department of Mathematics, Sungkyunkwan University, Suwon
440-746, Korea. \& School of Mathematics,
    Korea Institute for Advanced Study,
     Seoul 130-722, Korea} \email{parkj@skku.edu}
\begin{abstract} {We extend a result of Patodi for closed Riemannian
manifolds to the context of closed contact manifolds by showing the
condition that a manifold is an $\eta$-Einstein Sasakian manifold is
spectrally determined. We also prove that the condition that a
Sasakian space form has constant $\phi$-sectional curvature $c$ is
spectrally determined.
\\MSC 2010: 58J50, 53C25.
\\Keywords: $\eta$-Einstein manifold, Sasaki manifold, spectral geometry.}\end{abstract}

\maketitle

\section{ Introduction}\label{Sect1}
The relationship between the spectrum of certain natural operators
of Laplace type and the underlying geometry of a Riemannian manifold
has been studied by many authors.
Let
$(M,g)$ be a compact Riemannian manifold.  Let $\Delta_p$ be the
Laplace-Beltrami operator acting on the space of smooth $p$ forms
over a compact $m$-dimensional Riemannian manifold $M$.
 Patodi \cite{refPa70} established the
following spectral characterization of space forms:
\begin{thm}\label{thm1}   Let $(M_i,g_i)$ be compact Riemannian manifolds
without boundary. Assume that
$\Pspec(\Delta_{p},M_1)=\Pspec(\Delta_p,M_2)$ for $p=0,1,2$.
Then:\begin{enumerate}\item The manifold $M_1$ has constant scalar
curvature $c$ if and only if the manifold $M_2$ has constant scalar
curvature $c$.
\item The manifold $M_1$ is Einstein if and only if the manifold $M_2$ is Einstein.
\item The manifold $M_1$ has constant sectional curvature $c$ if and only if the manifold $M_2$ has
constant sectional curvature $c$.
\end{enumerate}\end{thm}
 Donnelly \cite{Do75} and Gilkey and Sacks
\cite{GiSa75} extended Theorem \ref{thm1} to the complex setting,
and the present author extended Theorem \ref{thm1} from the context
of closed Riemannian manifolds to the context of compact Riemannnian
manifolds with boundaries \cite{P1}. See also related work
\cite{P2}.

A contact metric manifold ${M}$ of dimension $m$ with contact form
$\eta$ and associated metric ${g}$ is called an {\it $\eta$-Einstein
manifold} if the Ricci tensor ${\rho}$ is given by
$${\rho} = \alpha {g} +
\beta \eta \otimes \eta\quad\text{for}\quad \alpha,\beta\in
C^\infty(M)\,.$$ Note that $\alpha$ and $\beta$ are constant if $M$
is a $\eta$-Einstein Sasakian manifold of dimension $\geq5$
\cite{Blair1}; this fails if $\text{dim} {M}=3$ \cite{JPS}. Also note that
the $\eta$-Einstein tangent sphere bundle of a Riemannian manifold
$M$ of radius $r$ equipped with the standard contact metric
structure has constant functions $\alpha$ and $\beta$ \cite{C1,
Park1} if $\text{dim} M\ge2$.

 The study of $\eta$-Einstein metrics is related to the Sasakian
Calabi problem \cite{BGK}. Tanno \cite{T} showed that Sasaki metric
on the unit tangent sphere bundle of any sphere $S^n$ is
$\eta$-Einstein and $D$-homothetic deformation of this metric
produces a homogeneous Einstein metric on $T_1{S^n}$. We refer to
\cite{D, HV, P, Z} for related work and some physical applications.
In this paper, we shall extend our study in the Riemannian setting
to the case of the contact geometry setting. The following is the
main result of this paper:
\begin{thm}\label{thm2}  Let $M_i = ({M_i},\eta_i,{g_i},\phi_i,\xi_i)$ be $m_i$-dimensional
compact Sasakian manifolds
without boundary with $m_i\ge5$. Assume
that $\Pspec(\Delta_{p}, M_1)=\Pspec(\Delta_{p}, M_2)$ for
$p=0,1,2$. Then:\begin{enumerate}
\item $m_1=m_2$ and $\operatorname{Vol}(M_1)=\operatorname{Vol}(M_2)$.
\item  $M_1 $ has constant scalar
curvature $c$ if and only if the manifold $M_2$ has constant scalar curvature
$c$.
\item  $M_1$ is $\eta$-Einstein
if and only if $M_2$ is
$\eta$-Einstein.
\item $M_1$ is Sasakian space form with constant $\phi$-sectional curvature $c$
if and only if $M_2$ is Sasakian
space form with constant $\phi$-sectional curvature
$c$.\end{enumerate}
\end{thm}

{\ The values $p=0,1,2$ are not particularly special. They are chosen for
 illustrative purposes only - there are other values which could be chosen -- see related
 work \cite{S, T82} for example. If $f\in C^\infty(M)$, let
 $$f[M]:=\int_Mf(x)\operatorname{dvol}(x)\,.$$
 The crucial point is that under the hypotheses of either Theorem~\ref{thm1} or of
 Theorem~\ref{thm2} that
 $$\{1[M],\tau[M],\tau^2[M],|\rho|^2[M],|R|^2[M]\}$$
 are spectrally determined.}\\

Here is a brief outline to the remainder of this paper. In Section
\ref{Sect2}, we review some facts concerning the Sasakian manifold.
In Section \ref{Sect3}, we review some previous results concerning
the heat trace asymptotics. In Section \ref{Sect4},
we complete the proof of Theorem \ref{thm2}.

\vskip.3cm The author would like to thank to Professors Gilkey and
Sekigawa for their helpful comments on the manuscript.

\section{Sasakian manifolds}\label{Sect2}

All manifolds in the present paper are assumed to be connected and
of class $C^{\infty}$.  We prepare some fundamental material about
Sasakian manifold. We refer to \cite{Blair1} for further details.
 A $(2n+1)$-dimensional manifold $M^{2n+1}$ is said to
be a {\em contact manifold} if it admits a global $1$-form $\eta$
such that $\eta\wedge(d\eta)^n\neq0$ everywhere. Given a contact
form $\eta$, we have a unique vector field $\xi$, the {\em
characteristic vector field,} satisfying $\eta(\xi)=1$ and
$d\eta(\xi,X)=0$ for any vector field $X$. It is well-known that
there exists a Riemannian metric $g$ and a $(1,1)$-tensor field
$\phi$ such that
\begin{equation}\label{2.1}
        \eta(X)=g(X,\xi),\quad d\eta(X,Y)=g(X,\phi Y),
        \quad \phi^2 X=-X+\eta(X)\xi
\end{equation}
where $X$ and $Y$ are vector fields on $M$. From (\ref{2.1}), it
follows that
\begin{equation}\label{2.2}
        \phi\xi=0,\quad \eta\circ\phi=0,
        \quad g(\phi X,\phi Y)=g(X,Y)-\eta(X)\eta(Y).
\end{equation}
A Riemannian manifold $M$ equipped with structure tensors
$(\eta,g,\phi,\xi)$ satisfying (\ref{2.1}) is said to be a {\em
contact metric manifold} (or {\em contact Riemannian manifold}) and
is denoted by $M=(M,\eta,g,\phi,\xi)$. \\

A normal contact metric manifold is called a {\it Sasakian}
manifold. Equivalently, an almost contact metric manifold ${M} =
({M},\eta,{g},\phi,\xi)$ is a Sasakian manifold if and only if the
following condition holds \cite{Blair1}:
 \begin{equation}\label{2.29}
({\nabla}_{{X}}\phi){Y} = {g}({X},{Y})\xi - \eta ({Y}) {X}.
 \end{equation}
On the other  hand, a contact metric manifold is called a
$K$-contact manifold if the characteristic vector field $\xi$ is a
Killing vector field. It is well-known that a Sasakian manifold is
necessarily a $K$-contact manifold.
 We also have the following formulas for a Sasakian manifold \cite{Blair1}.
\begin{equation}\label{2.30}
\nabla_X \xi=-\phi X,\\
\end{equation}
 \begin{equation}
 \begin{split}\label{2.31}
  R(X,Y)\xi= \eta(Y)X-\eta(X)Y,
 \end{split}
\end{equation}
 \begin{equation}\label{2.32}
\rho(\xi, \xi) = 2n,
\end{equation}
where $R$ and $\rho$ are the curvature tensor and Ricci tensor of
$M$, respectively.
\begin{defn}
{{A contact metric manifold ${M} = ({M},\eta,{g},\phi,\xi)$ is said
to be {\it $\eta$-Einstein} if the Ricci tensor $\rho$ of  $M$ is of
the form
$$\rho =\alpha g +\beta \eta \otimes \eta$$ for smooth functions
$\alpha$ and $\beta$ on $M$.}}
\end{defn}
\noindent On the other hand, it is known that any 3-dimensional
Sasakian manifold is $\eta$-Einstein. We may easily check that, for
an $\eta$-Einstein Sasakian manifold, the functions $\alpha$ and
$\beta$ are both constant if
 $\text{dim} M \geq 5$ \cite{Blair1}. {{A
Sasakian manifold $M$ is called a {\em Sasakian space form} if $M$
has constant $\phi$-sectional curvature. It is known that the
curvature tensor \cite{YanoKon} of a $2n + 1\geq 5$-dimensional
Sasakian space form with constant {\em $\phi$-sectional curvature}
is given by
\begin{equation}\label{5.1}
\begin{split}
R(X,Y,Z,W)&=g(R(X,Y)Z,W)\\
&=\frac{c+3}{4}\{g(Y,Z)g(X,W)-g(X,Z)g(Y,W)\large\}\\
&\;\;\;\;+\frac{c-1}{4}\{g(X,\phi Z)g(\phi Y,W)-g(Y,\phi Z)g(\phi X,
W)+2g(X,\phi Y)g(\phi Z,W)\}\\
&\;\;\;\;+\frac{c-1}{4}\{\eta(X)\eta(Z)g(Y,W)-\eta(Y)\eta(Z)g(X,W)+g(X,Z)\eta(Y)\eta(W)\\
&\;\;\;\;-g(Y,Z)\eta(X)\eta(W)\}.
\end{split}
\end{equation}
for any vector fields $X, Y, Z, W$ on $M$.
 Let $M$ be a
$(2n+1)$-dimensional Sasakian space form. We set $m = 2n + 1$. Then
from \eqref{5.1}, we see that the Ricci tensor $\rho$ of $M$ is
given by
\begin{equation}\label{505}
\rho = \frac{1}{4}\{(m+1)c + 3m -5\} g -
\frac{m+1}{4}(c-1)\eta\otimes\eta,
\end{equation}
 and hence, $M$ is an
$\eta$-Einstein manifold.

{{We now define the tensor fields $S_{\alpha, \beta}$  and $T_c$ of
$M$ respectively by
\begin{equation}\label{33}
S_{\alpha, \beta} (X, Y)= \rho(X, Y) - (\alpha g(X, Y) + \beta
\eta(X)\eta(Y)), \quad \text{and}
\end{equation}
\begin{equation}\label{34}
\begin{split}
T_c(X,Y,Z,W) = &R(X,Y,Z,W)
-\large\{\frac{c+3}{4}\{g(Y,Z)g(X,W)-g(X,Z)g(Y,W)\}\\
&\;\;+\frac{c-1}{4}\{g(X,\phi Z)g(\phi Y,W)-g(Y,\phi Z)g(\phi X,
W)+2g(X,\phi Y)g(\phi Z,W)\}\\
&\;\;+\frac{c-1}{4}\{\eta(X)\eta(Z)g(Y,W)-\eta(Y)\eta(Z)g(X,W)+g(X,Z)\eta(Y)\eta(W)\\
&\;\;-g(Y,Z)\eta(X)\eta(W)\}\large\}.
\end{split}
\end{equation}}}
for vector fields  $X, Y, Z, W$ on $M$, where $\alpha, \beta$ are
some smooth functions on $M$ and $c$ is a constant.\\

 Let $\{e_i\}$ be an orthonormal basis of
$T_pM$ at any point $p\in M$. In the sequel, we shall adopt the
following notational convention:
\begin{equation}
\begin{gathered}
    R_{ijkl}=g(R(e_i,e_j)e_k,e_l), \\
    \rho_{ij}=\rho(e_i,e_j),
    \quad
\phi_{ij}=g(\phi e_i, e_j), \quad\\
 \quad
 \nabla_i\phi_{jk}=g((\nabla_{e_i}\phi)e_j,e_k),
 \nabla_i\eta_{j}=g((\nabla_{e_i}\xi),e_j)
\end{gathered}
\end{equation}
and so on, where the Latin indices run over the range $1,2,\cdots,m
= 2n + 1$. We adopt the Einstein summation convention for the
repeated indices. From \eqref{2.29} to \eqref{2.32}, we may rewrite
as follows:
\begin{equation}\label{7}
\begin{split}
&\nabla_i \phi_{jk}=g_{ij}\eta_{k}-\eta_{j}g_{ik}, \\
&\nabla_i \eta_{j} =-\phi_{ij}, \\
&R_{ijkl}\eta_k = \eta_jg_{il} - \eta_i g_{jl},\\
&\rho_{ij}\eta_i \eta_j = 2n.
\end{split}
\end{equation}

\noindent From the definition of the tensor field $S_{\alpha,
\beta}$ and \eqref{7}, by direct calculation, we have
\begin{equation}\label{44}
|S_{\alpha, \beta}|^2 = |\rho|^2 - 2\alpha \tau + \gamma,
\end{equation}
where $\gamma = m \alpha^2 + 2\alpha\beta + \beta^2 - 2(m-1)\beta$.
Further, we see that $M$ is $\eta$-Einstein with the coefficient
functions $\alpha$ and $\beta$ in the defining equation if and only
if $S_{\alpha, \beta} = 0$ and $\alpha = \frac{\tau}{m-1} - 1$,
$\beta = m - \frac{\tau}{m -1}$ hold. In this case, we note that
$\alpha$, $\beta$ and $\tau$ are all constant if $\text{dim}M \geq
5$.

Next, we prepare the following Lemma to calculate
 the square norm $|T_c |^2$ of the tensor field $T_c$ on $M$.

\begin{lem}\label{lem2} On Sasakian manifold, we have
$$R_{ijkl}\{\phi_{ki}\phi_{jl}-\phi_{kj}\phi_{il}+2\phi_{ji}\phi_{kl}\}=6\tau-6(m-1)^2.$$
\end{lem}
\noindent{\bf{Proof }}
\noindent First, we get
\begin{equation}\label{4.1}
R_{ijkl}\phi_{ki}\phi_{jl}=\frac{1}{2}(R_{ijkl}-R_{kjil})\phi_{ki}\phi_{jl}=\frac{1}{2}(R_{ijkl}+R_{jkil})\phi_{ki}\phi_{jl} =-\frac{1}{2}R_{kijl}\phi_{ki}\phi_{jl}.\\
\end{equation}
Similarly, we obtain
\begin{equation}\label{4.2}
\begin{split}
-R_{ijkl}\phi_{kj}\phi_{il}&=-\frac{1}{2}R_{jkil}\phi_{jk}\phi_{il},\\
R_{ijkl}\phi_{ji}\phi_{kl}&=-R_{ijkl}\phi_{ij}\phi_{kl}.
\end{split}
\end{equation}
From \eqref{4.1} and  \eqref{4.2}, we have
\begin{equation}\label{4.3}
R_{ijkl}\{\phi_{ki}\phi_{jl}-\phi_{kj}\phi_{il}+2\phi_{ji}\phi_{kl}\}
= - 3 R_{ijkl}\phi_{ij}\phi_{kl}.
\end{equation}
On the other hand, from \eqref{7}, we get
\begin{equation*}
\nabla_l\nabla_i \phi_{jk}=g_{ij}\nabla_l
\eta_{k}-g_{ik}\nabla_l\eta_j=-g_{ij}\phi_{lk}+g_{ik}\phi_{lj},
\end{equation*}
and hence
\begin{equation}\label{8}
\nabla_l\nabla_i \phi_{jk}-\nabla_i\nabla_l
\phi_{jk}=-g_{ij}\phi_{lk}+g_{ik}\phi_{lj}+g_{lj}\phi_{ik}-g_{lk}\phi_{ij}.
\end{equation}
\noindent Applying the Ricci identity to (\ref{8}), and then taking
sum by setting $i=k$ in the resulting equality, we get
\begin{equation}\label{61}
-R_{lija}\phi_{ai}-\rho_{la}\phi_{ja} = (m-2)\phi_{lj}.
\end{equation}
Transvecting $\phi_{lj}$ to \eqref{61}, and taking account of
\eqref{7}, we have
\begin{equation}\label{29}
-R_{lija}\phi_{ai}\phi_{lj}-\rho_{la}\phi_{ja}\phi_{lj}=(m-2)\phi_{lj}\phi_{lj}=(m-1)(m-2),
\end{equation}
and hence
\begin{equation}\label{98}
-R_{ilaj}\phi_{ai}\phi_{jl} = -\rho_{la}(g_{la} - \eta_{l}\eta_{a})+
(m-1)(m-2)=-\tau + (m-1)^2.
\end{equation}
Thus, from \eqref{4.1} and \eqref{98}, we have
\begin{equation}\label{97}
\frac{1}{2}R_{jkil}\phi_{jk}\phi_{il} = -\tau + (m-1)^2.
\end{equation}
Therefore,  from \eqref{4.3} and \eqref{97}, we have
\begin{equation}\label{4.13}
R_{ijkl}\{\phi_{ki}\phi_{jl}-\phi_{kj}\phi_{il}+2\phi_{ji}\phi_{kl}\}
= 6\tau - 6(m-1)^2.
\end{equation}
This completes the proof of Lemma \ref{lem2}. \hfill$\square$\medskip\\

\section{Heat trace asymptotics}\label{Sect3}

Let $M$ be a compact Riemannian manifold of real dimension $m$
without boundary, and let $D$ be a operator of Laplace type on the
space of smooth sections to a smooth vector bundle over $M$. Let
$e^{-tD}$ be the fundamental solution of the heat equation. This
operator is of trace class and as $t\downarrow0$ there is a complete
asymptotic expansion with locally computable coefficients in the
form:
$$\Tr _{L^2}e^{-tD}\sim\textstyle\sum_{n\ge0}t^{(n-m)/2}a_n(D).$$

To study the heat trace coefficients $a_n(D)$, we introduce a bit of
additional notation. There is a canonically defined connection
$\nabla=\nabla(D)$ and a canonically defined endomorphism $E=E(D)$
so that
$$D=-(\Tr(\nabla^2)+E).$$

\noindent Let indices $i$, $j$, $k$ range from $1$ to $m$ and index
a local orthonormal frame $\{e_1,...,e_m\}$ for $TM$. Let $\Omega$
be the curvature of $\nabla$, let $\tau:=R_{ijji}$ be the scalar
curvature, let $\rho_{ij}:=R_{ikkj}$ be the Ricci tensor. Let `;'
denote multiple covariant differentiation. We refer to \cite{BrGi}
for the proof of the following result:

\begin{thm}\label{thm2.1} Let $D$ be an operator of Laplace type on the space of sections
$C^\infty(V)$ to a vector bundle $V$ over a compact manifold $M$.
Let $I$ be the identity endomorphism of $V$. We have:
\begin{enumerate}
\item $a_{0}(D)=(4\pi)^{-m/2}\int_M\Tr \{I\}$.
\smallbreak\item $a_{2}(D{})=(4\pi)^{-m/2}\frac16\int_M\Tr\{6E+\tau
I\}.$ \smallbreak\item
$a_{4}(D{})=(4\pi)^{-m/2}\frac1{360}\int_M\Tr\{60E_{;kk}+60\tau
E+180E^{2}+30\Omega^{2}$\smallbreak$
\;\;\;\;\;\;\;\;\;\;\;\;\;+(12\tau_{;kk} +
5\tau^{2}-2|\rho|^{2}+2|R|^{2})I\}.$
\end{enumerate}
\end{thm}
\noindent Theorem \ref{thm2.1} plays an important role in the proof
of Theorem \ref{thm2}. We refer to \cite{Gi03} for further
details.

\section{Proof of Theorem 2}\label{Sect4}

Let $M = ({M},\eta,{g},\phi,\xi)$ be a $2n + 1 \geq 5$-dimensional
compact Sasakian manifold without boundary, and set $m = 2 n +1$.
From Theorem \ref{thm2.1}, for $D=\Delta_p$ ($p = 0, 1, 2$), we have
\begin{equation}\label{4.41}
\textstyle\Tr _{L^2}(e^{-t\Delta_{0}})=(4\pi
t)^{-m/2}\{\vol(M)+O(t)\}\, \quad\text{and also}\nonumber
\end{equation}
\begin{equation}\label{4.42}
a_2(\Delta_0, M) = \frac16(4\pi)^{-m/2}\int_{M} \tau,
\end{equation}
\begin{equation}\label{4.43}
a_2(\Delta_1, M) = \frac16(4\pi)^{-m/2}\int_{M} {(m - 6)\tau},
\end{equation}
The work of Patodi \cite{refPa70} shows that there exist universal
constants so:
\begin{equation}\label{4.44}
a_4(\Delta_p, M) = (4\pi)^{-m/2}\int_M \{c^1_{m,p} \tau^2
+c_{m,p}^2|\rho|^2+c_{m,p}^3|R|^2+c_{m,p}^4\tau_{;ii}\}\,.
\end{equation}
 $p= 0, 1, 2$.

Now, we shall prove Theorem \ref{thm2}. Let $M_i =
({M_i},\eta_i,{g_i},\phi_i,\xi_i)$ be $m_i$-dimensional compact
Sasakian manifolds without boundary of $m_i \geq 5$ ($i = 1, 2$).
Assume that $\Pspec(\Delta_{p}, M_1)=\Pspec(\Delta_{p}, M_2)$ for
$p=0,1,2$. We denote by $R_i$, $\rho_i$ and ${\tau_i}$ the curvature
tensor, the Ricci tensor and the scalar curvature of $M_i$ ($i=1,
2$), respectively. Then, from $a_0(\Delta_0, M_i)$ in Theorem
\ref{thm2.1} (1), we have
\begin{equation}\label{41}
m_1 = m_2 \quad\text{and}\quad \vol(M_1)=\vol(M_2).
\end{equation}
We then establish assertion (2) by computing:
\begin{equation}\label{3.6}
\begin{split}
\tau_1&=(4\pi)^{m/2}\vol(M_1)^{-1}\{ma_2(\Delta_{0},M_1)-a_2(\Delta_{1},M_1)\}\\
&=(4\pi)^{m/2}\vol(M_2)^{-1}\{ma_2(\Delta_{0},M_2)-a_2(\Delta_{1},M_2)\}\\
&=\tau_2.
\end{split}
\end{equation}
The assertion (2) is nothing but a special case of the Theorem
\ref{thm1} (1).

Next, suppose that $M_1$ is an $\eta$-Einstein manifold with the
coefficient functions $\alpha_1$ and $\beta_1$ in the defining
equation. Since $m_1 \geqq 5$, it follows that $\alpha_1$ and
$\beta_1$ are constant and hence, the scalar curvature $\tau_1$ of
$M_1$ is also constant given by $\tau_1 = m\alpha_1 + \beta_1$.
Thus, from assertion (2), it follows that the scalar curvature
$\tau_2$ of $M_2$ is also constant and $\tau_{1}=\tau_{2}$. Since
$\vol(M_1)=\vol(M_2)$, the integrals of $\tau^2$ are equal. Since
$\tau_{;ii}=0$, from \eqref{4.44},
we have
\begin{equation}\label{jhp-eqn4}
\textstyle\int_{M_1}(c_{m,p}^2|\rho_1|^2+ c_{m,p}^3|R_1|^2)
=\textstyle\int_{M_2}(c_{m,p}^2|\rho_2|^2+ c_{m,p}^3|R_2|^2)
\end{equation}
for $p=1,2$; these two equations are independent \cite{refPa70}.
Consequently
\begin{equation}\textstyle\int_{M_1}|\rho_1|^2=\int_{M_2}|\rho_2|^2\quad\text{and}\quad
\int_{M_1}|R_1|^2=\int_{M_2}|R_2|^2.\label{eqn443}\end{equation}

\noindent Thus, from \eqref{44} we have
\begin{equation}\label{4.47}
0=\int_{M_1} |S^1_{\alpha_1, \beta_1}|^2=\int_{M_1}
|\rho_1|^2-2\alpha_1\tau_1 + \gamma_1 = \int_{M_2}
|\rho_2|^2-2\alpha_1\tau_1 + \gamma_1, \,
\end{equation}
 where $\gamma_1 = m \alpha{_1}^2 + 2\alpha_1\beta_1 + \beta_1^2 -
2(m-1)\beta_1$. Here, we may note that
\begin{equation}\label{4.48}
\begin{split}
&\alpha_1 = \frac{\tau_1}{m-1} - 1 =\frac{\tau_2}{m-1} - 1, \\
&\beta_1= m - \frac{\tau_1}{m -1} = m - \frac{\tau_2}{m -1}.
\end{split}
\end{equation}
We here set
\begin{equation}\label{4.50}
\begin{split}
S_{\alpha, \beta}^2 = \rho_2
 - (\alpha_2 g_2 + \beta_2 \eta_2 \otimes \eta_2 ),\quad
\text{where}\quad\alpha_2 = \frac{\tau_2}{m-1} - 1, \beta_2 = m -
\frac{\tau_2}{m -1}.
\end{split}
\end{equation}
 Then, we have
\begin{equation}\label{4.51}
\int_{M_2} |S^2_{\alpha_2, \beta_2}|^2=\int_{M_2}
|\rho_2|^2-2\alpha_2\tau_2 + \gamma_2, \,
\end{equation}
where $\gamma_2 = m \alpha{_2}^2 + 2\alpha_2\beta_2 + \beta_2^2 -
2(m-1)\beta_2.$\\

\noindent From \eqref{4.48} and \eqref{4.50}, we get
\begin{equation}\label{4.52}
\alpha_1 = \alpha_2,\quad  \beta_1 = \beta_2 \quad\text{and
hence}\quad \gamma_1 = \gamma_2.
\end{equation}
Therefore, from \eqref{4.47}, \eqref{4.51} and \eqref{4.52}, we have
$0=\int_{M_2} |S^2_{\alpha_2, \beta_2}|^2$, and therefore, $M_2$ is
an $\eta$-Einstein manifold with the same coefficients in the
defining equation. This completes the proof of Theorem \ref{thm2}
(3).\\


Lastly, suppose that $M_1$ is a Sasakian space form with constant
$\phi$-sectional curvature $c$. Then, from \eqref{505}, we see that
$M_1$ is an $\eta$-Einstein manifold with constant coefficients
$\alpha_1 = \frac{1}{4}\{(m+1)c + 3m -5\}$ and $\beta_1 = -
\frac{m+1}{4}(c-1)$ in the defining equation. Thus, it follows that
the scalar curvature $\tau_1$ is given by $$\tau_1 =
\frac{m-1}{4}\{(m+1)c + 3m -1\}.$$ Thus, by the assertion (3) and
hypothesis that $\Pspec(\Delta_{p}, M_1)=\Pspec(\Delta_{p}, M_2)$
(p=0, 1, 2), we see that $M_2$ is an $\eta$-Einstein manifold with
the constant coefficients $\alpha_2$ and $\beta_2$ in the defining
equation such that $\alpha_2 = \alpha_1$, $\beta_2 = \beta_1$, and
hence $\tau_2 = \tau_1$.
 We denote by $T_{c}^{2}$ the tensor field defined
by \eqref{34} of the Sasakian manifold $M_2$. We set
\begin{equation}\label{3.2}
(T_c^2)_{ijkl} := R_{ijkl} - K_{ijkl},
\end{equation}
where
\begin{equation*}
\begin{split}
K_{ijkl}&=\frac{c+3}{4}\{g_{jk}g_{il}-g_{ik}g_{jl}\}
+\frac{c-1}{4}\{\phi_{ki}\phi_{jl}-\phi_{kj}\phi_{il}+2\phi_{ji}\phi_{kl}\}\\
&\;\;\;\;+\frac{c-1}{4}\{\eta_i \eta_k g_{jl}-\eta_j \eta_k
g_{il}+g_{ik}\eta_j \eta_l-g_{jk}\eta_i \eta_l\}.
\end{split}
\end{equation*}
Then, by direct calculation, we have 

\begin{equation}\label{3.3}
\begin{split}
|K|^2=&\frac{(c+3)^2}{16} \{g_{jk}g_{il}-g_{ik}g_{jl}\}^2
+\frac{(c-1)^2}{16}(\phi_{ki}\phi_{jl}-\phi_{kj}\phi_{il}+2\phi_{ji}\phi_{kl})^2\\
&+\frac{(c-1)^2}{16}\eta_i \eta_k g_{jl}-\eta_j \eta_k g_{il}+
g_{ik}\eta_j \eta_l -g_{jk}\eta_i
\eta_l)^2\\
&+ \frac{(c-1)(c+3)}{8}
(g_{jk}g_{il}- g_{ik}g_{jl})(\phi_{ki}\phi_{jl}-\phi_{kj}\phi_{il}+2\phi_{ji}\phi_{kl})\\
&+ \frac{(c-1)(c+3)}{8} (g_{jk}g_{il}-g_{ik}g_{jl})(\eta_i \eta_k
g_{jl}-\eta_j \eta_k g_{il} + g_{ik}\eta_j \eta_l - g_{jk}\eta_i
\eta_l)\\
&+ 0 \\&=\frac{m-1}{2}\{(m+1)c^2 + 3m -1\}.
\end{split}
\end{equation}


Next, by taking account of Lemma \ref{lem2}, we get

\begin{equation}\label{77}
\begin{split}
R_{ijkl}K_{ijkl}=& \frac{c+3}{4}
R_{ijkl}(\delta_{jk}\delta_{il}-\delta_{ik}\delta_{jl})\\
&+\frac{c-1}{4}R_{ijkl}(\phi_{ki}\phi_{jl}-\phi_{kj}\phi_{il}+2\phi_{ji}\phi_{kl})\\
&+\frac{c-1}{4} R_{ijkl}(\eta_i \eta_k \delta_{jl}-\eta_j \eta_k
\delta_{il}+\delta_{ik}\eta_j \eta_l-\delta_{jk}\eta_i \eta_l)\\ &=
2c\tau - \frac12 (m-1)(3m -1)(c-1).
\end{split}
\end{equation}
Then, from \eqref{3.2}, \eqref{3.3} and \eqref{77}, we have
\begin{equation}\label{78}
|T_c^2|^2 = |R_2|^2 - 4c\tau_2 + d,
\end{equation}
where $d := \frac{m-1}{2}(m+1)c^2 + (m-1)(3m-1)c
-\frac{1}{2}(m-1)(3m -1)$. On the other hand, since $M_1$ is
$m$-dimensional Sasakian space form with constant $\phi$-sectional
curvature c, we have
$$|R_1|^2 = \frac{m-1}{2}\{(m+1) c^2 + 3m -1\},$$
and further
\begin{equation}\label{79}
0 =  |T_c^1|^2 = |R_1|^2- 4c\tau_1 + d.
\end{equation}
Then we use \eqref{eqn443}, \eqref{78}, taking account of \eqref{79}
and $\tau_1 = \tau_2$, we have
\begin{equation*}
0= \int_{M_1} |T_c^1|^2 =\int_{M_1} |R_1|^2- 4c\tau_1 + d =
\int_{M_2} |R_2|^2- 4c \tau_2 + d = \int_{M_2} |T_c^2|^2,
\end{equation*}
and hence, $T_c^2 = 0$ on $M_2$. Therefore, we see that $M_2$ is
also an $m$-dimensional Sasakian space form with constant
$\phi$-sectional curvature $c$.
This completes the proof of Theorem \ref{thm2} (4). \hfill$\square$\medskip\\

\section*{Acknowledgement}

\noindent This work was supported by the National Research
Foundation of Korea (NRF) grant funded by the Korea government
(MEST) (2011-0012987).


\end{document}